\documentclass{article}
\usepackage{amsmath}%
\usepackage{amssymb,latexsym,amsthm}
\usepackage{url}
\usepackage{geometry}
\usepackage[all]{xy}
\xyoption{all}

\newcommand {\iso}         {\cong}

\let\wt=\widetilde
\let\lra=\longrightarrow
\let\l=\lambda
\let\ol=\overline
\let\fr=\partial
\let\d=\delta
\let\e=\varepsilon
\let\phi=\varphi
\let\sub=\subseteq

\let\bs=\backslash
\let\inc=\hookrightarrow
\let\ck=\check
\let\mpty=\varnothing

\newcommand{\pfaff}{\mathrm{Pfaff}}

\newcommand{\E}{\mathcal{E}}

\newcommand{\dist}{\mathrm{dist}}

\let\Om=\Omega
\newcommand{\x}{\mathbf{x}}
\newcommand{\y}{\mathbf{y}}

\newcommand{\Int}{\mathrm{int}}
\newcommand{\im}{\mathop{\mathrm{im}}}

\newcommand{\R}{\ensuremath{\mathbb{R}}}
\newcommand{\RP}{\ensuremath{\mathbb{RP}}}
\newcommand{\N}{\ensuremath{\mathbb{N}}}

\renewcommand{\S} {\ensuremath{\mathcal{S}}} 
\renewcommand{\P} {\ensuremath{\mathcal{P}}}
\newcommand{\U} {\ensuremath{\mathcal{U}}}
\newcommand{\X} {\ensuremath{\mathcal{X}}}
\newcommand{\Y} {\ensuremath{\mathcal{Y}}}

\newcommand{\junk}[1]{}

\newcommand{\bbm}{b^{\mbox{\tiny{BM}}}}
\newcommand{\Hbm}{H^{\mbox{\tiny{BM}}}}

\newtheorem{theorem}{Theorem}[section]
\newtheorem{lemma}[theorem]{Lemma}
\newtheorem{corollary}[theorem]{Corollary}
\newtheorem{definition}[theorem]{Definition}
\newtheorem{proposition}[theorem]{Proposition}
\newtheorem*{main}{Main Result}
\theoremstyle{definition}
\newtheorem{example}[theorem]{Example}
\newtheorem{rem}[theorem]{Remark}
\newtheorem*{notations}{Notations}
\newtheorem*{ack}{Acknowledgements}

\begin{document}
\title
{Topological complexity %
of the relative closure \\ of a semi-Pfaffian couple}
\author{Thierry Zell}


\maketitle

\begin{abstract}
Gabrielov introduced the notion of relative closure of a Pfaffian
couple as an alternative construction of the o-minimal structure
generated by Khovanskii's Pfaffian functions. In this paper, we use the
notion of format (or complexity) of a Pfaffian couple to derive
explicit upper-bounds for the homology of its relative closure.
We consider both the singular and the Borel-Moore homology theories.

\end{abstract}

\section*{Introduction}

Pfaffian functions form a class of real-analytic functions with
finiteness properties similar to that of polynomials (see
$\mathsection1.1$). They were introduced by Khovanskii~\cite{kh91} who
proved for them an analogue of the theorem of B\'ezout: a system of
$n$ Pfaffian functions in $n$ variables can only have finitely many
isolated solutions.  In~\cite{wilkie99}, Wilkie proved that the
structure $\S_\pfaff$ generated by Pfaffian functions is {\em
  o-minimal}, thus confirming the intuition that sets defined using
such well-behaved functions must have tame topological properties.
(We refer the reader to~\cite{coste:pisa} and~\cite{vdd98} for more on
o-minimal structures and their topological properties.)


Pfaffian functions can be endowed with a notion of {\em complexity}
(known as {\em format}\/), a tuple of integers which can be used to
give an explicit upper-bound in Khovanskii's theorem, and, more
generally, to study quantitative aspects of the sets in $\S_\pfaff.$
Many estimates appear in the litterature, especially in the case {\em
semi-Pfaffian sets}, which are the sets defined by quantifier-free
Pfaffian formulas (Definition~\ref{df:sP}). A non-exhaustive list of
results about semi-Pfaffian sets 
would include the complexity of the frontier and
closure~\cite{ga:frcl} and of weak stratifications~\cite{gv:strat}, as
well as bounds on the sum of the Betti numbers~\cite{z99}.


In order to extend the notion of format to any definable set from
$\S_\pfaff$, Gabrielov introduced in~\cite{ga:rc} the notion of {\em
relative closure} $(X,Y)_0$ of a {\em semi-Pfaffian couple} $(X,Y)$
(see~$\mathsection$\ref{sub:rc}).  For the present introduction, it
suffices to say that a relative closure is definable set of
$\S_\pfaff$ constructed from the Hausdorff limits of two semi-Pfaffian
families $X$ and $Y$ depending on one parameter $\l$. The main result
in~\cite{ga:rc} is that any set in $\S_\pfaff$ is a finite union of
such relative closures.


%
Gabrielov's construction suggests a natural way to extend the
definition of format from semi-Pfaffian sets to general elements of
$\S_\pfaff$. Since the fibers $X_\l$ and $Y_\l$ of a semi-Pfaffian
couple $(X,Y)$ are semi-Pfaffian, we let   the format of the couple
$(X,Y)$ be the tuple which is the component-wise maximum of the
formats of $X_\l$ and $Y_\l$ (for some fixed small $\l>0$).
This definition leads to new quantitative
results, such as upper-bounds on the number of connected
components~\cite{gz} of a relative closure, and on the higher Betti
numbers~\cite{z03} under the assumption $Y=\mpty$.  In this paper, we
conclude this study of the Betti numbers of relative closures by
dealing with the case where $Y$ is not empty. 
Our results can be summarized in the following statement.

\begin{main}
Let $H_*$ and $\Hbm_*$ denote the homology groups associated
respectively to the singular and the Borel-Moore theories.
Then, for any semi-Pfaffian couple $(X,Y)$, 
the rank of the groups 
\begin{equation*}
H_k((X,Y)_0) \text{ and } \Hbm_k((X,Y)_0)
\end{equation*}
admit an upper-bound that is an explicit function of $k$ and of the
format of the couple $(X,Y)$.  In particular, the format of the
semi-Pfaffian sets $X$ and $Y$ in the parameter variable $\l$ {\em
does not appear} in these estimates.
\end{main}

We leave the detailed definitions and specific estimates until later
sections.  The Borel-Moore case (Theorem~\ref{th:RC-bm}) is a
reduction to the case $Y=\mpty$.
The singular case (Theorem~\ref{th:rcbetti}) is more involved: it
features a reduction to a definable Hausdorff limit of a family that
is {\em not} semi-Pfaffian. We then use an ad-hoc spectral sequence
argument to estimate the Betti numbers in that case.

The paper is organized as follows: Section~1 introduces the Pfaffian
structure and Gabrielov's construction of it via relative closures. It
also presents all the spectral sequence machinery and its corollaries
that appeared first in~\cite{gvz} and~\cite{z03}, and which will be
used in our proofs. Section~2 is devoted to the Borel-Moore estimates,
and Section~3 to the singular case.

\begin{notations}
For any $\x\in\R^n$, $|\x|$ denotes the euclidean norm of $\x$.
For any set $X\sub \R^n$, $\ol{X}$ denotes the closure of $X$  
in the euclidean topology, and $\fr X=\ol{X}\bs X$ denotes the
frontier of $X$. We also denote by $b_k(X)$ the rank of the homology
group $H_k(X)$, and let $b(X)=\sum_k b_k(X)$ be the sum of the Betti
numbers of $X$.
\end{notations}

\begin{ack}
The author is indebted to Andrei Gabrielov for help with
Lemma~\ref{lem:betticomp}, and to Nicolai Vorobjov for useful comments on
an earlier version of this paper.
\end{ack}

\section{Preliminaries}

In this section, we discuss Pfaffian functions and related notions:
semi-Pfaffian sets, the o-minimal structure $\S_\pfaff$ generated by
Pfaffian functions, and the description of $\S_\pfaff$ by relative
closures and limit sets. To each of these constructions, we can
associate a notion of complexity that we will call {\em format.} The
reader can find more details on Pfaffian sets and complexity results
in the survey~\cite{gv:survey}.

\subsection{Pfaffian functions}

Let $\U \sub \R^n$ be an open domain. The following definition is due
to Khovanskii~\cite{kh91}.

\begin{definition}
\label{df:chain}
Let $\x=(x_1, \ldots, x_n)$ and let $(f_1(\x), \ldots, f_\ell(\x))$ be
a sequence of analytic functions in $\U.$ This sequence is called a
{\em Pfaffian chain} if the functions $f_i$ are solution on $\U$ of a
triangular differential system of the form;
\begin{equation}\label{eq:chain}
df_i(\x)=\sum_{j=1}^n P_{i,j}(\x, f_1(\x), \ldots, f_i(\x))dx_j;
\end{equation}
where the functions $P_{i,j}$ are polynomials in $\x$ and $(f_1, \ldots,
f_i).$
\end{definition}

\begin{definition}
\label{df:Pfunc}
Let $(f_1, \ldots, f_\ell)$ be a fixed Pfaffian chain on a domain
$\U.$ The function $q$ is a {\em Pfaffian function expressible in the
chain} $(f_1, \ldots, f_\ell)$ if there exists a polynomial $Q$ such
that for all $\x \in \U,$
\begin{equation}\label{eq:Q}
q(\x)=Q(\x,f_1(\x), \ldots, f_\ell(\x)).
\end{equation}
In general, a function $q: \U \to \R$ is called {\em
Pfaffian} if it is expressible in some Pfaffian chain $(f_1, \ldots,
f_\ell)$ defined on $\U.$
\end{definition}


If $(f_1, \ldots, f_\ell)$ is a Pfaffian chain, we call
$\ell$ its {\em length}, and we let its {\em degree} $\alpha$ be the
maximum of the degrees of the polynomials $P_{i,j}$ appearing
in~\eqref{eq:chain}.  If $q$ is as in~\eqref{eq:Q}, the degree $\beta$
of the polynomial $Q$ is called the {\em degree} of $q$ in the chain
$(f_1, \ldots, f_\ell).$ 

\begin{definition}
\label{df:qformat}
For $q$ as above, the tuple $(n,\ell,\alpha,\beta)$ is called the {\em
format} of $q.$
\end{definition}

Pfaffian functions form a large class that contains, among other
things, real elementary functions and Liouvillian functions
(see~\cite{kh91}).  For more on the practical complexity of Pfaffian
functions, the papers~\cite{ga:cex,gv:strat,gv:survey} contain
examples of elementary functions, their various Pfaffian
representations, and the associated formats.

\subsection{Semi-Pfaffian sets}
We fix $(f_1, \ldots, f_\ell)$ a Pfaffian chain  defined on a domain
$\U,$ which we  assume to be of the form
\begin{equation}\label{eq:U}
\U=\{ \x \in \R^n \mid  g_1(\x)>0,\ldots, g_k(\x)>0\};
\end{equation}
where $g_1, \ldots, g_k$ are Pfaffian functions that are expressible
in the chain $(f_1, \ldots, f_\ell)$.

\begin{definition}
\label{df:qf}
Let $\P=\{p_1, \ldots, p_s\}$ be a set of Pfaffian functions
expressible in the chain $(f_1, \ldots, f_\ell)$.
A formula $\Phi$ is called a {\em quantifier-free formula on $\P$} if
it is derived from atoms of the form $p_i \star 0$ for $1 \leq i
\leq s$ and $\star \in \{=,<,>\}$, using conjunctions, disjunctions
and negations.
\end{definition}

\begin{definition}
\label{df:sP}
A subset $X \sub \R^n$ is called a {\em semi-Pfaffian set} if there
exists a quantifier-free Pfaffian formula $\Phi$ whose atoms are
Pfaffian functions expressible in some chain $(f_1, \ldots, f_\ell)$
defined on a domain $\U \sub \R^n$ of the form~\eqref{eq:U}
such that $X=\{\x \in \U \mid \Phi(\x)\}$.

The semi-Pfaffian set $X$ is called {\em restricted} if and only if
$\ol{X} \sub \U.$ 
\end{definition}

\begin{rem}
If $X$ is a restricted semi-Pfaffian set, its
closure is also semi-Pfaffian~\cite{ga:frcl}, a fact which is {\em
not} true for non-restricted sets~\cite{ga:cex}. This fact will be
used in Section~2.2.
\end{rem}

\begin{definition}
\label{df:format}
Let  $(f_1, \ldots, f_\ell)$ be a fixed Pfaffian chain and
$\P=\{p_1, \ldots, p_s\}$ be a collection of $s$ Pfaffian functions
expressible in that chain.
If the format
of each $p_i$ is bounded by  $(n,\ell,\alpha,\beta),$
then the {\em format} of any quantifier-free formula on $\P$, as well
as the
format of the corresponding semi-Pfaffian set,
is $(n,\ell, \alpha, \beta, s).$ 
\end{definition}


Khovanskii's estimate on the number of solutions of a system of
Pfaffian equations~\cite{kh91} allows, using arguments from Morse
theory, to bound the Betti numbers of any set which is defined as the common
zeros of a family of Pfaffian functions~\cite{kh91,z99}. This is used
to derive the following bound for semi-Pfaffian sets.

\begin{theorem}\label{th:qf}
Let $X$ be {\em any} semi-Pfaffian set defined by a quantifier-free
formula of format $(n,\ell, \alpha, \beta, s).$ The sum of the Betti
numbers of $X$ admits a bound of the form
\begin{equation}\label{eq:qf}
b(X) \leq 2^{\ell(\ell-1)/2} s^{2n} O(n (\alpha+ \beta))^{n+\ell};
\end{equation}
where the constant coming from the $O$ notation depends only on the
definable domain $\U$.
\end{theorem}
The estimate~\eqref{eq:qf} follows from the deformation techniques
used in~\cite{gv:qf} in the algebraic setting, by applying the bound
(for so-called $\P$-closed formulas) appearing in~\cite{z99}.  Note
that Theorem~\ref{th:qf} does not require to make assumptions either on the
topology of $X$ or on the formula defining $X$.

\subsection{Relative closure and limit sets}
\label{sub:rc}

Projections of semi-Pfaffian sets may not always be
semi-Pfaffian~\cite{ga:cex,osgood}, but Wilkie showed
in~\cite{wilkie99} that Pfaffian functions nonetheless generate an
o-minimal structure.
(See also~\cite{KaMcI,lr98,ls,s99} for related results.)
We  denote by $\S_\pfaff$ this o-minimal structure, and we call
(general) {\em Pfaffian set} any set which is definable in $\S_\pfaff$.
We refer the reader to~\cite{coste:pisa} and~\cite{vdd98}
for a  detailed account of the basic properties of o-minimal
structures.

In~\cite{ga:rc}, Gabrielov introduced the notions of relative closure
and limit sets to obtain a description of the structure $\S_\pfaff$
which allows to extend the notion of formats to all definable sets,
even those which are not definable by quantifier-free formulas.


To obtain all definable sets in $\R^n$, we need to consider 
semi-Pfaffian sets defined in a
domain $\U \sub \R^n \times \R_+$.  Without loss of generality, we
will assume that these sets are {\em bounded} (see Remark~\ref{rem:RPn}).
We write $\x=(x_1, \ldots ,x_n)$ for the coordinates in $\R^n$ and
$\lambda$ for the last coordinate (which we think of as a parameter).
If $X$ is  a subset of $\U$ and $\lambda>0,$ the {\em fiber} $X_\lambda$ is
defined by
\begin{equation*}
X_\lambda=\{\x \mid (\x,\lambda) \in X \} \sub \R^n;
\end{equation*}
and we consider $X$ as the family of its fibers $X_\lambda.$
We let $X_+=X \cap \{\lambda>0\}$ and denote by $\ck{X}$ the Hausdorff
limit of the family $\ol{X_\lambda}$ as $\lambda$ goes to zero;
\begin{equation}\label{eq:ck}
\ck{X}=\{\x \in \R^n \mid (\x,0)
\in \ol{X_+}\}.
\end{equation}

\begin{definition}
\label{df:family}
Let $X$ be a semi-Pfaffian subset of $\U$. The set $X$ constitutes a
{\em semi-Pfaffian family} if for any $\e>0,$ the set $X \cap
\{\lambda>\e\}$ is restricted. (Recall that by definition, the set $X
\cap \{\lambda>\e\}$ is restricted if its topological closure is
contained in $\U$.)
\end{definition}

\begin{definition}
\label{df:couple}
Let $X$ and $Y$ be semi-Pfaffian families in $\U$ defined in a common
chain $(f_1, \ldots, f_\ell).$ They form a {\em semi-Pfaffian couple}
if and only if, for all $\l>0$, we have $\ol{Y_\l}=Y_\l$ and
$\fr(X_\l)\sub Y_\l.$
\end{definition}


\begin{definition}
\label{df:formatcouple}
  The {\em format} $(n,\ell, \alpha, \beta, s)$ of a semi-Pfaffian
  family $X$ is the format of the fiber $X_\lambda$ for a small
  $\lambda>0.$ Then, the {\em format} of the couple $(X,Y)$ is the
  component-wise maximum of the format of the families $X$ and $Y.$
  \footnote{Note that the format of $X$ as a semi-Pfaffian {\em set}
    is different from its format as a semi-Pfaffian {\em family},
    since $X$ has $n+1$ variables as a set but only $n$ as a family.}
\end{definition}

\begin{definition}
Let $(X,Y)$ be a semi-Pfaffian couple in $\U.$ We define the {\em
relative closure} of $(X,Y)$ at $\lambda=0$ by 
\begin{equation}
\label{df:rc}
(X,Y)_0=\ck{X} \setminus \ck{Y} \sub \ck{\U};
\end{equation}
where $\ck{X}$, $\ck{Y}$ and $\ck{\U}$ denote the Hausdorff limits of
the respective fibers as in~\eqref{eq:ck}.
\end{definition}

\begin{rem}
\label{rem:X0}
The restrictions on semi-Pfaffian couples (Definition~\ref{df:couple})
imply that for $(X,\mpty)$ to be a couple, we must have
$\fr(X_\l)=\mpty$ for all $\l>0,$ and since $X$ is bounded, $X_\l$ must be
{\em compact}. We will denote by $X_0$ the relative closure
$(X,\mpty)_0$. In
that case, $X_0$ is simply the Hausdorff limit of the family of compacts
$X_\l$ when $\l$ goes to zero. 
\end{rem}

\begin{definition}
Let $\Om \sub \R^n$ be an open, semi-Pfaffian domain.  A {\em limit
set} $Z\sub\Om$ is a set of the form $Z=(X_1,Y_1)_0 \cup \cdots \cup
(X_k,Y_k)_0,$ where $(X_i,Y_i)$ are semi-Pfaffian couples respectively
defined in domains $\U_i \sub \R^n \times \R_+,$ such that
$\ck{\U_i}=\ol{\Om}$ for $1 \leq i
\leq k.$ If the formats of the couples $(X_i,Y_i)$ is bounded
component-wise by  $(n,\ell, \alpha, \beta, s)$ we say that the
{\em format} of the limit set is $(n,\ell, \alpha, \beta, s,k)$
\end{definition}

The main result of~\cite{ga:rc} is that limit sets are exactly the
definable sets in $\S_\pfaff$.
Moreover, the notion of format for limit sets makes the structure {\em
effective}, to the extent that when performing Boolean operations on
limit sets, the resulting formats can be explicitly bounded in terms
of the formats of the original sets.

\begin{rem}\label{rem:RPn}
When defining semi-Pfaffian couples, we assume, as in~\cite{ga:rc},
that the semi-Pfaffian families $X$ and $Y$ are bounded. This
restriction allows us to avoid a separate treatment of infinity: we
can see $\R^n$ as embedded in $\RP^n,$ in which case any set we
consider can be subdivided into pieces that are relatively compact in
their own charts.
\end{rem}

\subsection{Betti numbers of sub-Pfaffian sets and of Hausdorff limits}
\label{sub:Pi}

The first upper-bounds for the Betti numbers of Pfaffian sets which
were not defined by quantifier-free formulas were obtained
in~\cite{gvz} (for sub-Pfaffian sets) and in~\cite{z03} (for Hausdorff
limits). These results will be the key to our estimates for relative
closures.


Both results follow from the {\em descent inequality} presented
below. First, let's recall the following definition.

\begin{definition}
Let $X$ and $Y$ be two topological spaces, and let 
$f: X \to Y$ be a continuous surjection. The map $f$ is called {\em
locally split} if for any $y\in Y$, there exists a continuous
section of $f$ defined in a neighborhood of $y$.
\end{definition}
The following lemma gives a wide class of locally split maps.

\begin{lemma}\label{lem:locsplit}
Let $A$ and $B$ be topological spaces. If $U \sub A\times B$ is an
open subset for the product topology, then the restriction to $U$ of the
standard projection $\Pi: A \times B\to A$  is locally split.
\end{lemma}

\begin{proof}
Let $a^*\in \Pi(U)$. There must be $b^*\in B$ such that
$(a^*,b^*)\in U$, and since $U$ is open for the product topology,
there must be neighborhoods $V$ of $a^*$ and $W$ of $b^*$ such that
$V\times W\sub U$. The map $s: V \to V\times W$ given by
$s(a)=(a,b^*)$ is a continuous section of $\Pi|_U$; since this
construction can be done
for any $a^*\in \Pi(U)$, the map $\Pi|_U$ is locally split.
\end{proof}

\begin{theorem}\label{th:ss}
Let $f: X \to Y$ be a continuous surjective map definable in an
o-minimal structure.  Let $W^p_f(X)$ be the $(p+1)$-fold
fibered product of $X;$
\begin{equation}\label{eq:Wp-def}
W^p_f(X)=\{(\x_0,\ldots, \x_p) \in X^{p+1} \mid f(\x_0)=\cdots=f(\x_p)\}.
\end{equation}
Suppose that $f$ is either {\em closed} or {\em locally split}.
Then, for all $k$, the following inequality holds:
\begin{equation}\label{eq:th0}
b_k(Y) \leq \sum_{p+q=k} b_q(W^p_f(X)).
\end{equation}
\end{theorem}

The inequality~\eqref{eq:th0} follows from the existence of a spectral
sequence $E^r_{p,q}$ which converges to the homology of $Y$, and such that 
$E^1_{p,q}\iso H_q(W^p_f(X))$. The sequence $E^r_{p,q}$ is sometimes
known as the homological descent spectral sequence. It seems it first
appeared in~\cite{deligne} in the case of {\em proper}
maps, and it has been rediscovered many times since.
The reader can find proofs in~\cite{gvz} for the closed case and
in~\cite{bz} or~\cite{dugger} for the locally split case.

\begin{example}
We must note  that Theorem~\ref{th:ss} does not hold without
some kind of assumption on $f$, as the following example shows.  Let
$X \sub \R^3$ be the curve sketched in Figure~1, and let $f$ be the vertical
projection.
\begin{figure}
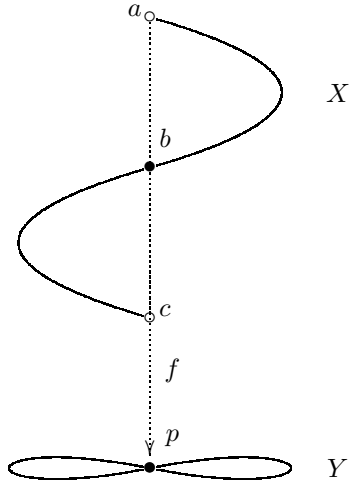
  
\[ \xy (0,30)*{\circ}="A"+(-2,1)*{a};
(0,10)*{\bullet}="B"+(2,4)*{b}; (0,-30)*{\bullet}="C"+(3,4)*{p};
(0,-10)*{\circ}="D"+(2,1)*{c}; (25,20)*{X}; (25,-30)*{Y}; (3,-17)*{f};
"A"; "B" **\crv{(35,20)}; "C"; "C" **\crv{(25,-25) & (25,-35)}; "C";
"C" **\crv{(-25,-25) & (-25,-35)}; "A"; "B" **\dir{.}; "B"; "D"
**\crv{(-35,0)}; "D"; "B" **\dir{.}; {\ar@{.>}"B";
(0,-30)*+{\bullet}};
\endxy  
\]
\caption{A map for which the inequality~\eqref{eq:th0} does not hold.}
\end{figure}
The curve $X$ is chosen so that $f$ is injective on $X$: the points
$a,b,$ and $c$ all project down to $p,$ but only $b$ is in $X.$ Since
$f$ is injective on $X,$ the set $W^1_f(X)$ is simply the diagonal in
$X^2,$ and since $X$ is contractible, we have
\begin{equation*}
b_0(W^1_f(X))+b_1(X)=1+0<b_1(Y)=2.
\end{equation*}
Thus, the inequality~\eqref{eq:th0} does not hold for $f$ when $k=1.$
\end{example}

\begin{definition}
A set $Y \sub \R^n$ is called {\em sub-Pfaffian} if there exists an
integer $r$ and a semi-Pfaffian set $X\sub \R^{n+r}$ such that
$Y=\Pi(X)$, where $\Pi$ is the standard projection $\R^{n+r} \to
\R^{n}$.
\end{definition}

When the inequality~\eqref{eq:th0} holds for such a projection $\Pi$,
we can apply Theorem~\ref{th:qf} and obtain a straightforward estimate
on the Betti numbers of the resulting sub-Pfaffian set.

\begin{corollary}\label{cor:Pi}
Let $X \sub \R^{n+r}$ be a semi-Pfaffian set of format $(n+r,\ell,
\alpha, \beta, s)$.  Denote by $\Pi$ the standard projection
$\R^{n+r} \to \R^{n}$ and let $Y=\Pi(X).$ Assume that the restriction
$\Pi|_X$ is either closed or locally split. Then, we have for all $k
\geq 1,$
\begin{equation}\label{eq:Pi}
b_{k-1}(Y) \leq 
2^{k\ell(k\ell-1)/2} s^{2(n+kr)} O((n+kr) (\alpha+ \beta))^{n+k(r+\ell)};
\end{equation}
\end{corollary}

\begin{proof}
The case where $\Pi|_X$ is closed is treated in~\cite{gvz}, and the proof is
identical in the locally split case.
\end{proof}


The relevance of Theorem~\ref{th:ss} to estimate the Betti
numbers of Hausdorff limits is not so obvious: 
it is the main result of~\cite{z03}, an estimate on the Betti numbers of
such a limit in terms of the Betti numbers of the associated expanded
diagonals (defined below).
\begin{definition}
\label{df:rho}
For any integer $p,$ we introduce the ``distance''
function $\rho_p$ on $(p+1)$-tuples $(\x_0, \ldots, \x_p)$ of points in
$\R^n$ by $\rho_0(\x_0)=0$, and for all $p\geq 1$,
\begin{equation}\label{eq:rho}
\rho_p(\x_0, \ldots,\x_p)=
\sum_{0 \leq i < j \leq p} |\x_i-\x_j|^2.
\end{equation}
For all $\e>0$ and
all integer $p \geq 0$, the {\em expanded $p$-th diagonal} of a set
$A\sub \R^n$ is then defined to be the subset of $(\R^n)^{p+1}$ given
by
\begin{equation}\label{eq:diag}
D^p(\e)=\{(\x_0, \ldots, \x_p) \in (A)^{p+1} \mid \rho_p(\x_0,
\ldots, \x_p)\leq \e\}.
\end{equation}
In particular, we have $D^0(\e)=A$ for all $\e>0.$
\end{definition}


The following result, which holds in any o-minimal structure, is
derived in~\cite{z03} from Theorem~\ref{th:ss}.

\begin{theorem}[{\cite[Theorem~1]{z03}}]
\label{th:haus}
Let $A \sub \R^{n+r}$ be a bounded set definable in some o-minimal
structure and let $A'$ be its projection to $\R^r.$ Suppose that the
fibers $A_a \sub \R^n$ are compact for all values of the parameter
$a\in A',$ and let $L$ be the Hausdorff limit of some sequence of
fibers $(A_{a_i}).$ Then, there exists $a\in A'$ and $\e>0$ such that
for any integer $k,$ we have
\begin{equation}\label{eq:bkL}
b_k(L) \leq \sum_{p+q=k} b_q(D^p_a(\e));
\end{equation}
where the set $D^p_a(\e)$ is the expanded $p$-th diagonal of the fiber
$A_a$.
\end{theorem}

We can apply Theorem~\ref{th:haus} 
to a special kind of relative closures. Indeed, we saw in
Remark~\ref{rem:X0} that if $(X,\mpty)$ was a semi-Pfaffian couple,
its relative closure $X_0=(X,\mpty)_0$ was simply the Hausdorff limit
as $\l$ goes to zero of the fibers $X_\l$ (which must be compact in
this case). Thus, Theorem~\ref{th:haus} shows that $b_k(X_0)$ can be
estimated in terms of the Betti numbers of the expanded diagonals,
which are semi-Pfaffian in this case. Applying Theorem~\ref{th:qf}, we
obtain the following explicit estimate.

\begin{corollary}[{\cite[Corollary~3]{z03}}]\label{cor:pfaff}
Let $X\sub \R^n \times \R_+$ be a semi-Pfaffian family with compact
fibers, and let $X_0=(X,\mpty)_0$ be the relative closure of $X.$ If the format of
$X_\l$ is bounded by $(n,\ell, \alpha, \beta, s),$ we have for any
$k\geq 1,$
\begin{equation}\label{eq:H-betti}
b_{k-1}(X_0) \leq  2^{k\ell(k\ell-1)/2} \, s^{2nk} \,
O(kn(\alpha+\beta))^{k(n+\ell)}.
\end{equation}
\end{corollary}

\section{Borel-Moore homology of relative closures}
\label{sec:borel}

In this section, we estimate the rank of the {\em Borel-Moore}
homology groups of the relative closure of a Pfaffian couple, in terms
of the format of the couple.  

\subsection{Borel-Moore homology in o-minimal structures}

In the o-minimal setting, the Borel-Moore homology of a locally
closed, definable set can be defined in the simple fashion described
below. We refer the reader to~\cite[$\mathsection$11.7]{bcr} for more
details.  

\begin{definition}
Let $S\sub \R^n$ be a set definable in some o-minimal structure.
If $S$ is compact, the {\em Borel-Moore homology} 
is simply $\Hbm_*(S)=H_*(S).$ If $S=A \bs B$ for some definable
compact sets $A$ and $B$ such that $B \sub A,$ the {\em Borel-Moore
homology} $\Hbm_*(S)$ is
\begin{equation}\label{eq:bmdef}
\Hbm_*(S)=H_*(A,B). 
\end{equation}
\end{definition}
Note that the Borel-Moore homology groups are {\em not} defined for
{\em all} definable subsets, but only for
{\em locally closed} sets, {\em i.e.} sets that can be written in the
form $U \cap F$ where $U$ is open and $F$ is closed.

We will denote by $\bbm_k(S)$ the rank of the group $\Hbm_k(S).$ The
Borel-Moore theory has an attractive property: the numbers
$\bbm_k$ are sub-additive.

\begin{proposition}[{\cite[Proposition~11.7.5]{bcr}}]\label{prop:subadd}
Let $S$ be a locally closed definable set and $T\sub S$ a closed
definable subset of $S.$  Then, there exists a long exact sequence 
\begin{equation*}
\cdots  \lra \Hbm_k(T) \lra \Hbm_k(S) \lra
\Hbm_k(S \bs T) \lra \Hbm_{k-1}(T) \lra\cdots
\end{equation*}
In particular, the following inequality hold for all integer $k;$
\begin{equation}\label{eq:subadd}
\bbm_k(S\bs T) \leq \bbm_k(S)+\bbm_{k-1}(T).
\end{equation}
\end{proposition}
This sub-additivity property is the key to bounding Borel-Moore ranks of
relative closures.

\subsection{Effective estimates in the Pfaffian structure}

Let us consider now a semi-Pfaffian couple $(X,Y)$.  For now, we will
assume that the fibers $X_\l$ and $Y_\l$ are compact for all $\l>0$,
and that their Hausdorff limits $\ck{X}$ and $\ck{Y}$ verify $\ck{Y}
\sub \ck{X}.$ Under these hypotheses, we have the following theorem.

\begin{theorem}\label{th:RC-bm}
Let $(X,Y)$ be a semi-Pfaffian couple as above. If the format of the
couple $(X,Y)$ is bounded by $(n,\ell, \alpha, \beta, s)$, then, for
any integer $k\geq 1,$ we have
\begin{equation}\label{eq:RC-bm}
\bbm_{k-1}((X,Y)_0) \leq
2^{k\ell(k\ell-1)/2} \, s^{2nk} \,
O(kn(\alpha+\beta))^{k(n+\ell)}.
\end{equation}
\end{theorem}

\begin{proof}
By definition, we have $(X,Y)_0=\ck{X}\bs \ck{Y}$. 
Since $\ck{X}$ and $\ck{Y}$ are compact sets such that $\ck{Y} \sub
\ck{X},$ Proposition~\ref{prop:subadd}  gives
\begin{equation*}%
\bbm_{k-1}((X,Y)_0)\leq %
\bbm_{k-1}(\ck{X})+\bbm_{k-2}(\ck{Y})
\end{equation*}
As the sets $\ck{X}$ and $\ck{Y}$ are compact, their Borel-Moore
homology coincides with the singular one, so that the previous
relation becomes
\begin{equation*}%
\bbm_{k-1}((X,Y)_0)\leq %
b_{k-1}(\ck{X})+b_{k-2}(\ck{Y}).
\end{equation*}
Since $\ck{X}$
(resp. $\ck{Y}$) is the
Hausdorff limit of the family of compact sets $X_\l$ (resp. $Y_\l$) 
when $\l$ goes to zero, the
ranks $b_{k-1}(\ck{X})$ and $b_{k-2}(\ck{Y})$ can be estimated using
Corollary~\ref{cor:pfaff}, yielding~\eqref{eq:RC-bm}.
\end{proof}

For a general semi-Pfaffian couple $(X,Y)$, two things can go wrong in
the previous argument:
the fibers $X_\l$ may not closed
\footnote{But remember that by definition of a semi-Pfaffian couple
(Definition~\ref{df:couple}),
the fibers $Y_\l$ have to be closed.}, and we may not have $\ck{Y} \sub
\ck{X}.$ 
If $\ck{Y} \not \sub \ck{X},$ we can simply consider the
couple $(X \cup Y, Y)$ which trivially verifies $(X \cup Y,
Y)_0=(X,Y)_0.$ The complexity of both couples is essentially the same,
and the inequality~\eqref{eq:RC-bm} still holds.


If the fibers $X_\l$ are not compact, a (weaker) bound can still be
established: since $X_\l$ is restricted, its closure $\ol{X_\l}$ is
also semi-Pfaffian, and its complexity can be estimated
using Theorem~1.1 of~\cite{ga:frcl}. Since taking the closure does not
change the Hausdorff limit $\ck{X},$ we can apply the above theorem to
the couple $(\ol{X},Y).$ However, the format of $\ol{X}$ involves
degrees that are doubly exponential in $n,$ so the bound on
$\bbm_{k-1}((X,Y)_0)$ is much worse than~\eqref{eq:RC-bm}.

\section{Singular homology of relative closures}
\label{sec:sing}

We will now establish a single-exponential upper-bound on the rank of
the singular homology groups of Pfaffian relative closures. Given a
semi-Pfaffian couple $(X,Y)$, we begin by constructing a family with
compact fibers $K$ which is definable in $\S_\pfaff$ (using one
universal quantifier) and whose Hausdorff limit has the same homology
groups as $(X,Y)_0$.

\subsection{Reduction to Hausdorff limits}
\label{sub:31}

Let $(X,Y)$ be a semi-Pfaffian couple and $(X,Y)_0$ be its relative
closure.  
If $\d: (0,1) \to \R_+$ is any function definable in $\S_\pfaff$, we let:
\begin{equation}\label{eq:K}
\begin{split}
K&=\{(\x,\l) \in X \cap \{\l>0\} \mid \dist(\x,Y_\l) \geq \d(\l)\}\\
&=\{(\x,\l) \in X \cap \{\l>0\} \mid \forall \y \in Y_\l, |\x-\y| \geq
\d(\l)\}. 
\end{split}
\end{equation}
The set $K$ is definable in $\S_\pfaff$; moreover, for any $\l>0$, the
set $K_\l$ is compact. Indeed, the fibers $X_\l$ are bounded by
assumption (see Remark~\ref{rem:RPn}) and $K_\l$ is a closed subset of
$X_\l$, so in order for $K_\l$ to be non-compact, it would have to
contain points in the frontier $\fr(X_\l)$. Equation~\eqref{eq:K}
shows that $K_\l \cap Y_\l=\mpty$, and as the definition of
semi-Pfaffian couple (Definition~\ref{df:couple}) requires that
$\fr(X_\l) \sub Y_\l$, the compactness of $K_\l$ follows.

\begin{proposition}\label{prop:K}
Let $(X,Y)$ be a semi-Pfaffian couple. Given a definable function
$\d(\l)$ as above and  $K$ defined as
in~\eqref{eq:K}, denote by $\d_0$ the limit of $\d(\l)$  
and by $K_0$ be the Hausdorff limit of the fibers
$K_\l$ (both limits taken when $\l$ goes to zero). Then, %
there exists $\d_1>0$ such that, 
for all $k\in\N$ and for all choice of $\d(\l)$ such that $\d_0\in(0,\d_1)$, 
we have
\begin{equation}\label{eq:bkK}
b_k((X,Y)_0)=b_k(K_0).
\end{equation}
\end{proposition}

\begin{proof}
Let $K(\d)$ be the definable subset given, for any $\d>0$, by
\begin{equation*}
K(\d)=\{\x \in \ck{X} \mid \dist(\x,\ck{Y}) \geq \d\}.
\end{equation*}
Since $(X,Y)_0=\{\x \in \ck{X} \mid \dist(\x,\ck{Y})>0\},$ the sets
$K(\d)$ are compact subset of $(X,Y)_0$ for all  $\d$ small enough.
The singular chain complexes $\{C_*(K(\d))\}$ form a directed system
for the inclusion maps. A singular chain of $(X,Y)_0$ has a support
which is contained in a set $K(\d)$ for a small enough value of $\d$;
thus, the direct limit of the system $\{C_*(K(\d))\}$ is
$C_*((X,Y)_0)$. By Theorem~4.1.7 of~\cite{spanier}, the homology and
direct limit functors commute, giving a similar equality for the singular
homology groups.
In particular, we must have
\begin{equation}\label{eq:limbk}
b_k((X,Y)_0)=\lim_{\d \to 0} b_k(K(\d)).
\end{equation}

Since the family $K(\d)$ is definable in $\S_\pfaff$, the generic
triviality theorem (see  Theorem~1.2 in Chapter~9 of~\cite{vdd98} or
Theorem~5.22 in~\cite{coste:pisa}) guarantees that we can find some real
number $\d_1>0$ such that the topological type of the sets $K(\d)$
is constant for $\d\in (0,\d_1).$ The relation~\eqref{eq:limbk}
implies that
\begin{equation*}
\forall \d\in (0,\d_1), \quad b_k((X,Y)_0)= b_k(K(\d)).
\end{equation*}
For the family $K$ defined in~\eqref{eq:K}, the Hausdorff limit $K_0$
of the fibers $K_\l$ as $\l$ goes to zero is $K_0=K(\d_0)$. Thus, any
definable function $\d(\l)$ that verifies $\d_0\in(0,\d_1)$ yields a
family $K$ such that~\eqref{eq:bkK} holds.
\end{proof}

Since $K_0$ is the Hausdorff limit of the definable family $K_\l$ when
$\l$ goes to zero, we can use Theorem~\ref{th:haus} to bound the Betti
numbers of $K_0$. Thus, there exists $\l>0$ and $\e>0$ such that
\begin{equation}\label{eq:K0}
b_k(K_0)\leq \sum_{p+q=k} b_q(D_\l^p(\e));
\end{equation}
where $D_\l^p(\e)$ denotes the
expanded diagonals of $K_\l$,
\begin{equation}\label{eq:DpK}
D_\lambda^p(\e)=\{(\x_0, \ldots, \x_p) \in (K_\lambda)^{p+1} \mid 
\rho_p(\x_0, \ldots, \x_p) \leq \e\}.
\end{equation}
The set $D_\lambda^p(\e)$ is defined by a Pfaffian formula which is
the conjunction of a quantifier free part
\begin{equation}\label{qfp}
(\x_0, \ldots, \x_p)\in (X_\l)^{p+1}
\wedge \rho_p(\x_0, \ldots, \x_p) \leq \e\;
\end{equation}
and a part using a single universal quantifier,
\begin{equation}\label{uqp}
\forall \y\in Y_\l, \quad \bigwedge_{i=0}^p |\x_i-\y| \geq \d(\l).
\end{equation}

\subsection{Complements and duality}
\label{sub:32}

Proposition~\ref{prop:K} and inequality~\eqref{eq:K0} reduce our
problem to bounding the Betti numbers of the sets $D_\l^p(\e)$, or
equivalently (via Alexander duality), of their complements.
The sets $(\U_\l)^{p+1}\bs D_\l^p(\e)$ are defined by existential
Pfaffian formulas; our next result estimates their Betti numbers using
the results from Section~\ref{sub:Pi}



\begin{proposition}\label{prop:betticomp}
Let $(X,Y)$ be a semi-Pfaffian family of format bounded by $(n,\ell,
\alpha, \beta, s)$ defined in a domain $\U.$ Let $p$ be some fixed
integer, $\l$ and $\e$ be positive real numbers, and let $D_\l^p(\e)$
be the set defined in~\eqref{eq:DpK}. For any integer $q$, the
Betti number $b_{q-1}((\U_\l)^{p+1}\bs D_\l^p(\e))$ is bounded by
%
\begin{equation}\label{eq:betticomp}
2^{[q(p+2)\ell]^2/2}\  [s(p+1)]^{O((p+q)n)}\ 
[n(p+q)(\alpha+\beta)]^{O((p+q)(n+\ell))}.
\end{equation}
\end{proposition}
%
%
%
%
%
\begin{proof}
In order to simplify notations, we let $\Om=(\U_\l)^{p+1},$
$D=D_\l^p(\e),$ $\d=\d(\l)$, $\X=(X_\l)^{p+1}$ and $\Y=Y_\l.$
With these notations, $D$ is the set of tuples $(\x_0, \ldots, \x_p)$
such that 
\begin{equation*}
(\x_0, \ldots, \x_p) \in \X\ \wedge\  \rho_p(\x_0, \ldots, \x_p)\leq
 \e\ \wedge\ 
 \forall \y \in \Y,\ \bigwedge_{i=0}^p |\x_i-\y_i| \geq \d.
\end{equation*}

Let $\Pi: \Om\times \Y \to \Om$ denote the projection on the first
factor, and let $\E=A\cup B \cup C$, where
\begin{align*}
A &=\bigcup_{i=0}^p\{(\x_0, \ldots, \x_p, \y)
\in \Om\times \Y \mid |\x_i-\y| <\d\};\\ %
B &=\{(\x_0, \ldots, \x_p, \y)
\in \Om\times \Y
\mid \rho_p(\x_0, \ldots, \x_p)>\e\}; \\
C&=\{(\x_0, \ldots, \x_p, \y)
\in \Om\times \Y \mid (\x_0, \ldots, \x_p)\not\in \X\};
\end{align*}
We have $\Om \bs D= \Pi(\E)$. 
In order to use the estimates on the Betti numbers of sub-Pfaffian
sets appearing in Corollary~\ref{cor:Pi}, we need 
to prove that
the restriction of $\Pi$ to $\E$ is locally split.
In order to do so, we introduce a shrinkage $\wt{C}$ of $C$ defined by
\begin{equation*}
\wt{C}=\{(\x_0, \ldots, \x_p, \y)
\in \Om\times \Y \mid (\x_0, \ldots, \x_p)\not\in \ol{\X}\}.
\end{equation*}
The point $(\x_0, \ldots, \x_p, \y)$ is in $C\bs
\wt{C}$ if and only if $(\x_0, \ldots, \x_p)\in \fr\X$, where 
\begin{equation*}
\fr\X=\bigcup_{i=0}^p \{(\x_0, \ldots, \x_p)\in \Om \mid \x_i\in \fr(X_\l)\}.
\end{equation*}
Since $(X,Y)$ is a semi-Pfaffian couple, we know that
$\fr(X_\l)\sub Y_\l$ (see Definition~\ref{df:couple}), so 
$(\x_0, \ldots, \x_p)\in \fr\X$ if and only if there exists $0 \leq i
\leq p$ such that $\dist(\x_i,Y_\l)=0$. In particular, this means that
$\Pi(C\bs \wt{C})=\fr\X \sub \Pi(A)$, giving
\begin{equation*}
\Pi(A\cup B \cup C)=\Pi(A\cup B \cup \wt{C}).
\end{equation*}
Let $\wt{\E}=A \cup B\cup\wt{C}$.
%
It is clear from the definition that $A,B$ and $\wt{C}$  are open
subsets of $\Om\times \Y$, so $\wt{\E}$ is open too, and according to
Lemma~\ref{lem:locsplit}, the restriction of $\Pi$ to $\wt{\E}$ is
locally split. But the restriction of
$\Pi$ to $\E$ must be locally split too: the local sections of
$\Pi|_{\wt{\E}}$ are local sections for $\Pi|_{\E}$, and  since
$\Pi(\E)=\Pi(\wt{\E})$, it is not necessary to check for the existence
of other sections. %
\footnote{For our argument, it is {\em not
    enough} to note that $\Pi|_{\wt{\E}}$ is locally split and to apply
  the spectral sequence in that case: indeed, both $\E$ and $\wt{\E}$
  are semi-Pfaffian, but the bound on the format of $\wt{\E}$ 
is much worse than the bound on the format of $\E$.}


%

If the format of $(X,Y)$ is bounded by $(n,\ell,
\alpha, \beta, s),$ the format of $\E$ is bounded by 
\begin{equation*}
(n(p+2),(p+2)\ell, \alpha, \max(2,\beta), (s+1)(p+2)).
\end{equation*}
Since $\E$ is semi-Pfaffian and $\Pi|_{\E}$ is locally split,
the Betti numbers of $\Pi(\E)=\Om \bs D$ are bounded according to
Corollary~\ref{cor:Pi}, and~\eqref{eq:betticomp} follows.
\end{proof}



In order to prove our main theorem, we will also need to relate the
Betti numbers of the set $(\U_\l)^{p+1} \bs D_\l^p(\e)$ to those of
the complement of $D_\l^p(\e)$ . This is achieved with the following
lemma.

\begin{lemma}\label{lem:betticomp}
Let $D$ and $\Om$ be subsets of $\R^N$ such that $\ol{D}\sub \Int(\Om).$
Then, for all $q$, we have 
\begin{equation*}
b_q(\R^N \bs D) \leq b_q(\Om \bs D).
\end{equation*}
\end{lemma}

\begin{proof}
To prove the result, it is enough to show that the map $k:\ H_q(\Om \bs X)
\to H_q(\R^N \bs X)$ induced by inclusion is surjective.
Let us consider the following commutative diagram, where the rows are
the exact sequences associated to the couples $(\R^N \bs D, \Om \bs
D)$ and $(\R^N,\Om)$ respectively, and the vertical arrows are induced
by the corresponding inclusions.
\begin{equation*}
\xymatrix{
 \cdots \ar[r]
&H_{q+1}(\R^N \bs D, \Om \bs D) \ar[d]_{\iso}^{i} \ar[r]^{\qquad\d}
& H_q(\Om \bs D) \ar[r]^{k} \ar[d]^{j}
& H_q(\R^N \bs D) \ar[r]^{\ell\quad} \ar[d]
& H_q(\R^N \bs D, \Om \bs D) \ar[d] \ar[r]^{\qquad \d} & \cdots \\
 \cdots \ar[r] 
& H_{q+1}(\R^N, \Om)  \ar[r]^{\fr}_{\iso}
& H_q(\Om) \ar[r]
& H_q(\R^N) \ar[r]
& H_q(\R^N, \Om) \ar[r] & \cdots 
}
\end{equation*}
Since $\ol{D} \sub \Int(\Om),$ the excision axiom 
asserts that the
inclusion $i: (\R^N \bs D, \Om \bs D) \inc (\R^N,\Om)$ is an isomorphism
on the homology level. Since $\R^N$ is contractible, the boundary maps
$\fr$ in the exact sequence of the couple $(\R^N,\Om)$ are
isomorphisms; thus, we obtain that the composition $\fr\circ i:
H_{q+1}(\R^N \bs D, \Om \bs D) \to H_q(\Om)$ is an isomorphism, and
since this map is equal to $j \circ\d,$ the map $\d$ must be
injective.

By exactness of the first row at $H_q(\R^N \bs D, \Om \bs D)$, we have
$\im\ell=\ker\d=0$ (since $\d$ injective), but by exactness at $H_q(\R^N
\bs D)$, we obtain $\ker \ell=H_q(\R^N \bs D)=\im k,$ and thus $k$ is
surjective. 
\end{proof}

\subsection{Betti numbers of a relative closure}


Our main result can now be obtained by combining the reduction of
Section~\ref{sub:31} and the estimates of Section~\ref{sub:32}.

\begin{theorem}\label{th:rcbetti}
Let $(X,Y)$ be a semi-Pfaffian family defined in a domain $\U\sub
\R^n\times \R_+$. 
If the  format of $(X,Y)$ is bounded by $(n,\ell,
\alpha, \beta, s)$, then for any integer $k\geq
1,$ the Betti number $b_{k}((X,Y)_0)$ is bounded by
\begin{equation}\label{eq:rcbetti}
2^{O(n^2k^4\ell^2)}\ [sk]^{O(kn)}\ [nk(\alpha+\beta)]^{O(k(n+\ell))};
\end{equation}
where the constants depend only on the domain $\U.$
\end{theorem}

\begin{proof}
We established from Proposition~\ref{prop:K} and Theorem~\ref{th:haus} 
the bound 
\begin{equation}
b_{k}((X,Y)_0)=b_{k}(K_0)\leq \sum_{p+q=k} b_q(D_\l^p(\e));
\end{equation}
for some suitable values $\l>0$ and $\e>0.$ 

Fix a value for $p.$ 
We denote by $N=(p+1)n$ the dimension of the ambient space containing
$D=D_\l^p(\e)$, and we use the notations introduced in
the proof of Proposition~\ref{prop:betticomp}. 
Since $D$ is compact, Alexander duality 
(Theorem~6.2.16 in~\cite{spanier}) gives $H^q(D)\iso
\tilde{H}_{N-q-1}(\R^N \bs D).$ Since $\Om$ is open and $D$ is closed,
Lemma~\ref{lem:betticomp} 
yields
\begin{equation*}
b_q(D) \leq b_{N-q-1}(\R^N \bs D)\leq b_{N-q-1}(\Om \bs D).
\end{equation*}
According to Proposition~\ref{prop:betticomp}, and since $p+q=k$, the
Betti number $b_{N-q-1}(\Om \bs D)$ is bounded by
\begin{equation}\label{eq:bqd}
2^{O(n^2p^4\ell^2)}\ [s(p+1)]^{O(kn)}\ [nk(\alpha+\beta)]^{O(k(n+\ell))}.
\end{equation}
Thus $b_q(D^p_\l(\e))$ is bounded by~\eqref{eq:bqd} too, and the sum 
$\sum_{p+q=k} b_q(D_\l^p(\e))$
is bounded by~\eqref{eq:rcbetti}, which proves the
theorem.
\end{proof}

\end{document}